\documentclass{article}
\usepackage{amssymb}
\usepackage{amsfonts}
\usepackage{amsmath}

\newtheorem{theorem}{Theorem}

\newtheorem{corollary}[theorem]{Corollary}

\newtheorem{lemma}[theorem]{Lemma}

\begin{document}

\title{Multitype branching processes evolving in i.i.d. random environment:
probability of survival for the critical case}
\author{Vatutin V.A.\thanks{%
Department of Discrete Mathematics, Steklov Mathematical Institute, 8,
Gubkin str., 119991, Moscow, Russia; e-mail: vatutin@mi.ras.ru}, Dyakonova
E.E.\thanks{%
Department of Discrete Mathematics, Steklov Mathematical Institute, 8,
Gubkin str., 119991, Moscow, Russia; e-mail: elena@mi.ras.ru}}
\date{}
\maketitle

\begin{abstract}
Using the annealed approach we investigate the asymptotic behavior of the
survival probability of a critical multitype branching process evolving in
i.i.d. random environment. We show under rather general assumptions on the
form of the offspring generating functions of particles that the probability
of survival up to generation $n$ of the process initiated at moment zero by
a single particle of type $i$ is equivalent to $c_{i}n^{-1/2}$ for large $n,$
where $c_{i}$ is a positive constant. Earlier such asymptotic representation
was known only for the case of the fractional linear offspring generating
functions.
\end{abstract}



\section{Introduction and main results}

\setcounter{equation}{0} \renewcommand{\theequation}{\thesection.%
\arabic{equation}} Branching processes in random environments (BPRE's) is an
important part of the theory of branching processes. This model was
introduced by Smith and Wilkinson \cite{SmWil} in 1969. Since then a lot of
papers appeared analyzing various properties of this and more general models
of BPRE's. As a result we now have a relatively complete description of the
basic properties of many models of BPRE's under the annealed and quenched
approaches. In particular, it became clear that the behavior of the
single-type BPRE's is mainly determined be the properties of the so-called
associated random walk constructed by the logarithms of the expected
population sizes of particles of the respective generations. For instance,
using the notion of the associated random walk it is easy to divide in a
natural way the set of all single-type BPRE's into the classes of
supercritical, critical and subcritical processes (see \cite{4h} for more
detail). Moreover, attracting known and proving new results for the ordinary
random walk on the real line it is possible to prove a wide range of
statements describing the long range behavior of the single-type BPRE's
under the annealed and quenched approaches.

The obtained results include theorems concerning the asymptotic behavior of
the survival probability of the processes up to a distant moment $n$,
describing the distribution of the population size of the BPRE at moments $%
nt,0\leq t\leq 1,$ given survival of the processes up to the moment $n$ (see
\cite{Af}, \cite{4h} -- \cite{ABKV2014}, \cite{BansVat}, \cite{BDKV10}, \cite%
{DyakGVa}, \cite{Kozlov}, \cite{VD} -- \cite{VK08}, \cite{Vat2012}, \cite%
{Vat2012b}, \cite{VW09}). In addition, limit theorems were proved which
describe the trajectories of the critical and subcritical BPRE's attaining a
high level (see \cite{Af2009} -- \cite{Af2014}), include a number of large
deviation type results (see \cite{BansBeres} --  \cite{BansBoe12}, \cite%
{BoeKerst13}, \cite{GrQuen2016}, \cite{HuaLiu12}, \cite{Koz06}, \cite{Koz10}%
, \cite{LiaLiu13}) and evaluate the rate of convergence to the limit of
supercritical processes (see \cite{HuaLiu14}, \cite{HuaLiu14b}). The reader
can find some other references and details in surveys \cite{VZ}, \cite%
{VDS2013} and \cite{Vat15}. Thus, we have now a relatively complete
description of the basic properties of the single-type BPRE's.

It is natural to try\ to prove analogues of the statements known for the
single-type case for the multitype BPRE's. This is, however, in no way
simple. The point is that the role of the associated random walk is played
in this case by the logarithms of the norms of products of certain random
matrices and many useful results known for ordinary random walks on the real
line have no analogues in terms of $p\times p$ random matrices with $p\geq 2.
$ \ Still, there are several papers dealing with multitype BPRE's. We
mention here only some of them: \cite{AthKar}, \cite{Dyak9}, \cite{Dyak},
\cite{VDEN16} -- \cite{Dayk13}, \cite{Kap}, \cite{Tanny},  \cite{Vat2010} --
\cite{Vat2010b},  \cite{VDS2013} -- \cite{VQuan2015}, \cite{Weis}.

Quite recently an important step in describing the asymptotic behavior of
the survival probability for the critical multitype BPRE under general
conditions has been made in \cite{LPP2016}. Using the ideas of papers \cite%
{GK}, \cite{Dyak9} and \cite{GPP2015} the authors of \cite{LPP2016} proved
that the probability of survival up to a distant moment $n$ has the order $%
n^{-1/2}$ for a wide class of multitype branching processes evolving in iid
random environment. Moreover, they have established that if the offspring
generating functions are fractional linear with probability 1 and the
process is initiated at time zero by a single particle of type $i,$ then the
probability of survival of the process up to moment $n$ behaves like $%
c_{i}n^{-1/2}$ as $n\rightarrow \infty $, where $c_{i}$ is a positive
constant. In the present paper we prove that the claims formulated in \cite%
{Dyak9} and \cite{LPP2016} for the survival probability of the multitype
BPRE with fractional linear offspring generating functions are valid for
more general classes of branching processes. Our prove is based on ideas contained in \cite{Dyak9}, \cite{Dyak}, \cite%
{GK}, \cite{GPP2015} and \cite{LPP2016}.

To formulate our main result we need some standard notation for
multidimensional vectors. In the sequence we usually make no difference in
notation for row and column vectors. It will be clear from the context which
form is selected.

Let $\mathbf{e}_{j},j=1,...,p,$ be the $p-$dimensional vector whose $j- $th
component is equal to $1$ and the others are zeros, $\mathbf{0}=(0,...,0) $
be the $p-$dimensional vector all whose components are zeros.

For $\mathbf{x}=(x_{1},...,x_{p})$ and $\mathbf{y}=(y_{1},...,y_{p})$ set%
\begin{equation*}
\left\vert \mathbf{x}\right\vert =\sum_{i=1}^{p}|x_{i}|,\quad (\mathbf{x},%
\mathbf{y})=\sum_{i=1}^{p}x_{i}y_{i}.
\end{equation*}%
For a $p\times p$ matrix $A=\left( A\left( i,j\right) \right) _{i,j=1}^{p}$,
all whose elements are nonnegative, set%
\begin{equation*}
|A|=\sum_{i,j=1}^{p}A\left( i,j\right) .
\end{equation*}%
Denote $\mathbb{N}_{0}=\{0,1,2,...\}$ and let $\mathbb{N}_{0}^{p}$ be the
set of all vectors $\mathbf{z}=(z_{1},...,z_{p})$ with non-negative
integer-valued coordinates. Let, further, $\mathbb{J}^{p}$ be the set of all
vectors
\begin{equation*}
\mathbf{s}=(s_{1},...,s_{p}),\quad ,0\leq s_{i}\leq 1,\,i=1,...,p.
\end{equation*}%
For $\mathbf{s}=(s_{1},...,s_{p})$ $\in \mathbb{J}^{p}$ and $\mathbf{z}%
=(z_{1},...,z_{p})\in \mathbb{N}_{0}^{p}$ set $\mathbf{s}^{\mathbf{z}%
}=\prod\limits_{i=1}^{p}s_{i}^{z_{i}}.$

Denote by $(\mathbb{F},\mathcal{B(}\mathbb{F}))$ the space of probability
measures $\Phi $ on $\mathbb{N}_{0}^{p}$ with $\sigma -$algebra $\mathcal{B(}%
\mathbb{F})$ of Borel sets endowed with the metric of total variation, and
let $(\mathbb{F}_{p},\mathcal{B(}\mathbb{F}_{p}))$ be the $p-$times product
of the space $(\mathbb{F},\mathcal{B(}\mathbb{F}))$ on itself.

Let $\mathbf{F=(}F^{(1)},...,F^{(p)})$ be a random vector (a tuple of random
measures) taking values in $(\mathbb{F}_{p},\mathcal{B(}\mathbb{F}_{p})),$
i.e., $\mathbf{F}$ is a measurable map of a probability space $(\Omega ,%
\mathcal{A},\mathbf{P}\mathcal{)}$ in $(\mathbb{F}_{p},\mathcal{B(}\mathbb{F}%
_{p})).$ \ An infinite sequence $\Pi =(\mathbf{F}_{0},\mathbf{F}_{1},\mathbf{%
F}_{2},...)$ of independent identically distributed copies of $\mathbf{F}$
is called a random environment and we say in this case that $\mathbf{F}$
generates $\Pi .$

Taking into account a one-to-one correspondence between probability measures
and generating functions we associate with $\mathbf{F=(}F^{(1)},...,F^{(p)})$
generating $\Pi $ a random $p-$dimensional vector $\mathbf{f}(\mathbf{s}%
)=(f^{(1)}\left( \mathbf{s}\right) ,...,f^{(p)}\left( \mathbf{s}\right) ),\,$
whose components are $p- $dimensional (random) generating functions $f^{(i)}$
corresponding to $F^{(i)},1\leq i\leq p:$
\begin{equation}
f^{(i)}(\mathbf{s}):=\sum_{\mathbf{z}\in \mathbb{N}_{0}^{p}}F^{(i)}(\{%
\mathbf{z}\})\mathbf{s}^{\mathbf{z}},\,\mathbf{s}\in \mathbb{J}^{p}.
\label{Defff}
\end{equation}%
In a similar way we associate with the component $\mathbf{F}%
_{n}=(F_{n}^{(1)},...,F_{n}^{(p)}),\,n\geq 0,$ of the random environment $%
\Pi $ a (random)vector of probability generating functions $\mathbf{f}_{n}(%
\mathbf{s})=(f_{n}^{(1)}\left( \mathbf{s}\right) ,...,f_{n}^{(p)}\left(
\mathbf{s}\right) ),$ whose components are (random) \ multi-variate
generating functions $f_{n}^{(i)}(\mathbf{s})$ corresponding to $%
F_{n}^{(i)},1\leq i\leq p,$
\begin{equation}
f_{n}^{(i)}(\mathbf{s}):=\sum_{\mathbf{z}\in \mathbb{N}_{0}^{p}}F_{n}^{(i)}(%
\{\mathbf{z}\})\mathbf{s}^{\mathbf{z}},\;\mathbf{s}\in \mathbb{J}^{p}.
\label{p.2}
\end{equation}

With this agreement in view we may formally rewrite (\ref{Defff}) and (\ref%
{p.2}) as

\begin{equation*}
f^{(i)}(\mathbf{s})=\mathbf{E}\left[ \mathbf{s}^{\mathbf{\xi }_{i}}|\mathbf{f%
}\right] =\mathbf{E}\left[ s_{1}^{\xi _{i1}}\cdot \cdot \cdot s_{p}^{\xi
_{ip}}|\mathbf{f}\right] ,\mathbf{s}\in \mathbb{J}^{p},
\end{equation*}%
and%
\begin{equation*}
f_{n}^{(i)}(\mathbf{s})=\mathbf{E}\left[ \mathbf{s}^{\mathbf{\xi }_{i}(n)}|%
\mathbf{f}_{n}\right] =\mathbf{E}\left[ s_{1}^{\xi _{i1}(n)}\cdot \cdot
\cdot s_{p}^{\xi _{ip}(n)}|\mathbf{f}_{n}\right] ,\mathbf{s}\in \mathbb{J}%
^{p},
\end{equation*}%
where, given $(F^{(1)},...,F^{(p)})=(\Phi ^{(1)},...,\Phi ^{(p)})$ the
vectors $\mathbf{\xi }_{i}=\left( \xi _{i1},...,\xi _{ip}\right) ,i=1,...,p,$
are distributed according to the measures $\Phi ^{(i)},i=1,...,p,$ and,
given $(F_{n}^{(1)},...,F_{n}^{(p)})=(\Phi _{n}^{(1)},...,\Phi _{n}^{(p)})$
the vectors $\mathbf{\xi }_{i}(n)=\left( \xi _{i1}(n),...,\xi
_{ip}(n)\right) ,i=1,...,p,$ are distributed according to the measures $\Phi
_{n}^{(i)},i=1,...,p,$ respectively.

In addition, we sometimes use $\Pi =(\mathbf{f}_{0},\mathbf{f}_{1},\mathbf{f}%
_{2},...)$ for $\Pi =(\mathbf{F}_{0},\mathbf{F}_{1},\mathbf{F}_{2},...)$.

Let $(\mathbf{\Phi }_{0},\mathbf{\Phi }_{1},\mathbf{\Phi }_{2},...),$ where
\begin{equation*}
\mathbf{\Phi }_{n}=(\Phi _{n}^{(1)},\Phi _{n}^{(2)},...,\Phi _{n}^{(p)})\in
\mathbb{F}_{p},\,n=0,1,2,...,
\end{equation*}%
be a realization of the random environment.

A sequence of random $p-$dimensional vectors $\mathbf{Z}%
(n)=(Z^{(1)}(n),...,Z^{(p)}(n)),$ $\;n=0,1,...,$ with non-negative
integer-valued coordinates is called a $p-$type Galton-Watson branching
process in random environment $\Pi ,$ if $\mathbf{Z}\left( 0\right) $ is
independent of $\Pi $ and for all $n\geq 0,\,\mathbf{z}=(z_{1},...,z_{p})\in
$ $\mathbb{N}_{0}^{p}$ and $\mathbf{\Phi }_{0},\mathbf{\Phi }_{1},...\in
\mathbb{F}_{p}$

\begin{equation}
\mathcal{L(}\mathbf{Z}\left( n+1\right) \;|\;\mathbf{Z}%
(n)=(z_{1},...,z_{p}),\Pi =(\mathbf{\Phi }_{0},\mathbf{\Phi }_{1},...))=%
\mathcal{L}\left( \sum_{i=1}^{p}\sum_{j=1}^{z_{i}}\mathbf{\xi }_{i}^{(j)}(n)|%
\mathbf{\Phi }_{n}\right) ,  \label{(1.1)}
\end{equation}%
where the tuple $\mathbf{\xi }_{i}^{(1)}(n),\mathbf{\xi }_{i}^{(2)}(n),...,%
\mathbf{\xi }_{i}^{(z_{i})}(n),i=1,...,p,$ consists of independent $p-$%
dimensional random vectors with non-negative integer-valued components.
Besides, for each $i\in \left\{ 1,...,p\right\} \,$the random vectors $%
\mathbf{\xi }_{i}^{(1)}(n),\mathbf{\xi }_{i}^{(2)}(n),...,\mathbf{\xi }%
_{i}^{(z_{i})}(n)$ are distributed according to the probability measure $%
\Phi _{n}^{(i)}$ $.$

Relation (\ref{(1.1)}) specifies a Galton-Watson branching process\ in
random environment in which $Z^{(i)}(n)$ is treated as the number of type $i$
particles in the $n-$th generation. Particles in the process evolve as
follows. If\textbf{\ }the state of the environment at time\textbf{\ }$n$ is $%
(\Phi _{n}^{(1)},...,\Phi _{n}^{(p)})\in \mathbb{F}_{p},$ then each
particles of type $i\in \left\{ 1,...,p\right\} $ existing at this moment in
the population produces offspring in accordance with the $p-$dimensional
probability measure $\Phi _{n}^{(i)}$ and independently of the reproduction
of other existing particles at time $n$ and of the prehistory of the
process. Thus, the component $Z^{(i)}(n+1)$ of the vector $\mathbf{Z}%
(n+1)=(Z^{(1)}(n+1),$ $...,Z^{(p)}(n+1))$ is equal to the number of type $i$
particles among all direct descendants of the particles of the $n-$th
generation.

In the sequel we assume that our BPRE is specified on a sufficiently reach
probability space $(\Omega ,\mathcal{A},\mathbf{P}\mathcal{)}$ and denote by
$\mathbf{E}\left[ \cdot \right] $ the expectation with respect to the
probability measure $\mathbf{P}$.

Besides, we use the symbols $\mathbb{E}_{\mathbf{f}}$ and $\mathbb{P}_{%
\mathbf{f}}$ for the expectations and probabilities given the offspring
generating function $\mathbf{f}=\mathbf{f}(\mathbf{s})=(f^{(1)}\left(
\mathbf{s}\right) ,...,f^{(p)}\left( \mathbf{s}\right) ),$ $\ \mathbf{s}\in
\mathbb{J}^{p},$ whose components are multivariate probability generating
functions. We also use the symbols $\mathbb{E}_{\mathbf{f}_{n}}$ and $%
\mathbb{P}_{\mathbf{f}_{n}}$ for the expectations and probabilities given
the environment $\mathbf{f}_{n}.$

Set $m_{ij}:=\mathbb{E}_{\mathbf{f}}\left[ \xi _{ij}\right] $ and $%
m_{ij}(n):=\mathbb{E}_{\mathbf{f}}\left[ \xi _{ij}(n)\right] ,i,j=1,...,p$
and introduce the corresponding mean matrices
\begin{eqnarray*}
M &=&M(\mathbf{f})=\left( m_{ij}\right) _{i,j=1}^{p}:=\left( \mathbb{E}_{%
\mathbf{f}}\left[ \xi _{ij}\right] \right) _{i,j=1}^{p}, \\
&& \\
M_{n} &=&M_{n}(\mathbf{f}_{n})=\left( m_{ij}(n)\right) _{i,j=1}^{p}:=\left(
\mathbb{E}_{\mathbf{f}_{n}}\left[ \xi _{ij}(n)\right] \right)
_{i,j=1}^{p},n\geq 0.
\end{eqnarray*}%
Along with random matrices $M_{n}$ we use, for $i\in \left\{ 1,...,p\right\}
$ the Hessian matrices%
\begin{eqnarray*}
B^{(i)} &=&B^{(i)}(\mathbf{f})=\left( B^{(i)}(k,l)\right)
_{k,l=1}^{p}:=\left( \mathbb{E}_{\mathbf{f}}\left[ \xi _{ij}(\xi
_{ik}-\delta _{kj})\right] \right) _{k,l=1}^{p}, \\
&& \\
B_{n}^{(i)} &=&B^{(i)}(\mathbf{f}_{n})=\left( B_{n}^{(i)}(k,l)\right)
_{k,l=1}^{p}:=\left( \mathbb{E}_{\mathbf{f}_{n}}\left[ \xi _{ij}(n)(\xi
_{ik}(n)-\delta _{kj})\right] \right) _{k,l=1}^{p}
\end{eqnarray*}%
and let%
\begin{eqnarray}
\mu &:&=\sum_{i=1}^{p}\left\vert B^{(i)}\right\vert ,\qquad \eta :=\frac{\mu
}{\left\vert M\right\vert ^{2}},  \label{DefEta} \\
&&  \notag \\
\mu _{n} &=&\sum_{i=1}^{p}\left\vert B_{n}^{(i)}\right\vert ,\qquad \eta
_{n}=\frac{\mu _{n}}{\left\vert M_{n}\right\vert ^{2}}.  \label{DefEtan1}
\end{eqnarray}

We define the cone
\begin{equation*}
\mathcal{C=}\left\{ \mathbf{x}=(x_{1},...,x_{p})\in \mathbb{R}^{p}:x_{i}\geq
0\text{ for any }i=1,...,p\right\} ,
\end{equation*}%
the sphere%
\begin{equation*}
\mathbb{S}^{p-1}=\left\{ \mathbf{x}:\mathbf{x}\in \mathbb{R}^{p},\left\vert
\mathbf{x}\right\vert =1\right\} ,
\end{equation*}%
and the space $\mathbb{X}=\mathcal{C}\cap \mathbb{S}^{p-1}$. To go furhter
we need to consider the general linear semigroup $S^{+}$ of $p\times p$
matrices all whose elements are non-negative and endow the semi-group with
the $L_{1}$ norm denoted also by $\left\vert \cdot \right\vert $. For $%
\mathbf{x}\in \mathbb{X}$ and $A\in S^{+}$ we specify the projective action
as%
\begin{equation*}
\mathbf{x}\cdot A:=\frac{\mathbf{x}A}{\left\vert \mathbf{x}A\right\vert }
\end{equation*}%
and define on the product space $\mathbb{X}\times S^{+}=\left\{ \left(
\mathbf{x},A\right) \right\} $ a function $\rho $ by setting
\begin{equation*}
\rho (\mathbf{x},A):=\ln \left\vert \mathbf{x}A\right\vert .
\end{equation*}%
For $A_{1},A_{2}\in S^{+}$ the function meets the so-called cocycle property%
\begin{equation*}
\rho (\mathbf{x},A_{1}A_{2})=\rho (\mathbf{x\cdot }A_{1},A_{2})+\rho (%
\mathbf{x},A_{1}).
\end{equation*}

The measure $\mathbf{P}$ introduced above for the multitype BPRE specifies
the respective probability measure on the Borel $\sigma $-algebra of the
semi-group $S^{+}$ which we also denote by $\mathbf{P}$, i.e., for a Borel
subset $\mathcal{B}\subseteq S^{+}$ we set%
\begin{equation*}
\mathbf{P}\left( M\in \mathcal{B}\right) =\mathbf{P}\left( \mathbf{f}:M=M(%
\mathbf{f})\in \mathcal{B}\right) .
\end{equation*}

Keeping in mind this agreement we introduce a number of assumptions to be
valid throughout the paper. These assumptions are taken from \cite{LPP2016}
and concern only the properties of the restriction of $\mathbf{P}$ to the
semi-group $S^{+}$.

\textbf{Condition H1}. There exists $\varepsilon _{0}>0$ such that $%
\int_{S^{+}}\left\vert M\right\vert ^{\varepsilon _{0}}\mathbf{P}\left(
dM\right) <\infty .$

\textbf{Condition H2}. (Strong irreducibility).\ The support of $\mathbf{P}$
in $S^{+}$ acts strongly irreducibly on $\mathbb{R}^{p},$ i.e. no proper
finite union of subspaces of $\mathbb{R}^{p}$ is invariant with respect to
all elements of the multiplicative semi-group generated by the support of $%
\mathbf{P}$.

\textbf{Condition H3}. There exists a real positive number $b>1$ such that
the elements of the matrix $M=M(\mathbf{f})=\left( m_{ij}\right)
_{i,j=1}^{p} $ meet the condition
\begin{equation*}
\frac{1}{b}\leq \frac{m_{ij}}{m_{kl}}\leq b
\end{equation*}%
for any $i,j,k,l\in \left\{ 1,...,p\right\} $.

It is know (see \cite{BL85}) that under Conditions H1-H3 there exists a
unique $\mathbf{P}$-invariant measure $\mathbf{v}$ on $\mathbb{X}$ such
that, for any continuous function $\varphi $ on $\mathbb{X}$,%
\begin{equation*}
\left( \mathbf{P}\ast \mathbf{v}\right) \left( \varphi \right)
:=\int_{S^{+}}\int_{\mathbb{X}}\varphi \left( \mathbf{x}\cdot M\right)
\mathbf{v}\left( d\mathbf{x}\right) \mathbf{P}\left( dM\right) =\int_{%
\mathbb{X}}\varphi \left( \mathbf{x}\right) \mathbf{v}\left( d\mathbf{x}%
\right) =:\mathbf{v}\left( \varphi \right) .
\end{equation*}

Moreover, if we denote by $R_{n}:=M_{0}M_{1}...M_{n-1}$ the right product of
random matrices $M_{k},k\geq 0,$ then (see \cite{FK60}), given%
\begin{equation}
\mathbf{E}\left[ \max (0,\ln \left\vert M\right\vert )\right] <\infty
\label{Fursten}
\end{equation}%
the sequence
\begin{equation*}
\frac{1}{n}\ln \left\vert R_{n}\right\vert ,n=0,1,2,...
\end{equation*}%
converges $\mathbf{P}$-a.s. to a limit
\begin{equation*}
\pi :=\lim_{n\rightarrow \infty }\frac{1}{n}\mathbf{E}\left[ \ln \left\vert
R_{n}\right\vert \right]
\end{equation*}%
called the upper Lyapunov exponent and, under Conditions H1-H3 (see \cite%
{BL85})
\begin{equation*}
\pi =\int_{S^{+}}\int_{\mathbb{X}}\rho \left( \mathbf{x},M\right) \mathbf{v}%
\left( d\mathbf{x}\right) \mathbf{P}\left( dM\right) .
\end{equation*}%
We introduce two more conditions.

\textbf{Condition H4}. The upper Lyapunov exponent of the distribution
generated by $\mathbf{P}$ on $S^{+}$ is equal to 0.

\textbf{Condition H5}. There exists $\delta >0$ such that

\begin{equation*}
\mathbf{P}\left( M\in S^{+}:\text{ for any }\mathbf{x}\in \mathbb{X},\,\ln
\left\vert \mathbf{x}M\right\vert \geq \delta \right) >0.
\end{equation*}

Now we are ready to formulate the main result of the paper.

\begin{theorem}
\label{T_main}Assume Conditions $H1-H5$. If
\begin{equation}
\sup_{\mathbf{x}\in \mathbb{X}}\mathbf{E}\left[ \frac{1}{\left\vert \mathbf{x%
}M\right\vert }\right] <\infty  \label{ExponFinite}
\end{equation}%
and, for $\eta $ given in (\ref{DefEta}) and some $\varepsilon >0$
\begin{equation}
\mathbf{E}\left[ \eta ^{1+\varepsilon }\right] <\infty ,
\label{SecondFinite}
\end{equation}%
then, for any $i\in \left\{ 1,...,p\right\} $ there exists a real number $%
\beta _{i}\in \left( 0,\infty \right) $ such that%
\begin{equation}
\lim_{n\rightarrow \infty }\sqrt{n}\mathbf{P}\left( \mathbf{Z}(n)\neq
\mathbf{0|Z(}0\mathbf{)=e}_{i}\right) =\beta _{i}.  \label{SurvivalProbab}
\end{equation}
\end{theorem}

This theorem refines the main results of \cite{LPP2016} where relation (\ref%
{SurvivalProbab}) was proved under the assumption that the offspring
generating functions are fractional linear and where it was shown for the
general case that if there exists a constant $C>0$ such that, for any $%
i,j,k\in \left\{ 1,...,p\right\} $%
\begin{equation*}
\mathbb{E}_{\mathbf{f}}\left[ \xi _{ik}(n)(\xi _{il}(n)-\delta _{kl})\right]
\leq C\mathbb{E}_{\mathbf{f}}\left[ \xi _{ik}\right] <\infty \text{ \qquad }%
\mathbf{P}\text{-a.s.}
\end{equation*}%
then, for any $i\in \left\{ 1,...,p\right\} $ there exist positive constants
\thinspace $c_{1}$ and $c_{2}$ such that%
\begin{equation*}
c_{1}\leq \sqrt{n}\mathbf{P}\left( \mathbf{Z}(n)\neq \mathbf{0|Z(}0\mathbf{%
)=e}_{i}\right) \leq c_{2}
\end{equation*}%
for all $n\geq 1.$

Theorem \ref{T_main} also complements the respective theorem from \cite%
{Vat2010b} where the asymptotic expression for the probability $\mathbf{P}%
\left( \mathbf{Z}(n)\neq \mathbf{0|Z(}0\mathbf{)=e}_{i}\right) $ is found
for the general form of the offspring distributions of particles of
different types under the following hypotheses:

$A1)$ the projection to $\mathbf{P}$ on $S^{+}$ is concentrated on a subset $%
\mathcal{B}(\mathbf{u})=\left\{ M\right\} \subset S^{+}$ of matrices having
a common nonrandom positive right eigenvector $\mathbf{u,}\left\vert \mathbf{%
u}\right\vert =1,$ corresponding to their maximal (in absolute value)
eigenvalues%
\begin{equation*}
M\mathbf{u}=\rho \mathbf{u.}
\end{equation*}

$A2)$ The sequence $\Sigma _{n}=\ln \rho _{0}+...+\ln \rho _{n-1}$ generated
by the sums of the logarithms of the maximal (in absolute value) eigenvalues
of the random matrices $M_{0},...,M_{n-1}$ satisfies the Doney-Spitzer
condition:%
\begin{equation*}
\lim_{n\rightarrow \infty }\mathbf{P}\left( \Sigma _{n}>0\right) =a\in
\left( 0,\infty \right) \text{.}
\end{equation*}%
In particular, the expectation $\mathbf{E}[\ln \vert M\vert] $ may
not exist meaning that condition (\ref{Fursten}) may not be valid.

Besides, Theorem \ref{T_main} complements the respective
result from \cite{Dyak} where the asymptotic expression for the probability $%
\mathbf{P}\left( \mathbf{Z}(n)\neq \mathbf{0|Z(}0\mathbf{)=e}_{i}\right) $
was obtained for the general form of the offspring distributions of
particles of different types under the following assumptions:

$A1^{\prime })$ the projection of $\mathbf{P}$ to $S^{+}$ is concentrated on
a subset $\mathcal{B}(\mathbf{v})=\left\{ M\right\} \subset S^{+}$ of
matrices having a common nonrandom positive left eigenvector $\mathbf{v,}%
\left\vert \mathbf{v}\right\vert =1,$ corresponding to their maximal (in
absolute value) eigenvalues, i.e.,%
\begin{equation*}
\mathbf{v}M_{n}=\rho _{n}\mathbf{v.}
\end{equation*}

$A2)$ The sequence $\Sigma _{n}$ satisfies the Doney-Spitzer condition.

\section{Auxiliary results}

Let $\left\{ \mathbf{X}_{n},n\geq 0\right\} $ be a homogeneous Markov chain
on $\mathbb{X}$ specified by the equalities%
\begin{equation*}
\mathbf{X}_{0}=\mathbf{x}\in \mathbb{X}\text{ and }\mathbf{X}_{n}=\mathbf{%
x\cdot }R_{n}.
\end{equation*}%
If the matrices $M_{n\text{ }}$ are i.i.d and subject to the probability law
$\mathbf{P}$ then for any $\mathbf{x}\in \mathbb{X}$ and any bounded Borel
function $\varphi $ the transition probabilities of this chain are given by
the formula
\begin{equation*}
Q\varphi (\mathbf{x})=\int_{S^{+}}\varphi \left( \mathbf{x}\cdot M\right)
\mathbf{P}\left( dM\right) .
\end{equation*}

For the subsequent arguments we need to consider a random process $\left\{
S_{n},n\geq 0\right\} $ defined for $\mathbf{x}\in \mathbb{X}$ and $a\in
\mathbb{R}$ by the relations
\begin{equation*}
S_{0}=S_{0}\left( \mathbf{x},a\right) :=a,\ S_{n}=S_{n}\left( \mathbf{x}%
,a\right) :=a+\ln \left\vert \mathbf{x}R_{n}\right\vert .
\end{equation*}%
According to the cocycle property%
\begin{equation*}
S_{n}\left( \mathbf{x},a\right) =a+\sum_{k=0}^{n-1}\rho (\mathbf{X}%
_{k},M_{k}).
\end{equation*}

The sequence of pairs $\left( \mathbf{X}_{n},S_{n}\right) ,n\geq 0,$ is a
Markov chain on $\mathbb{X}\times \mathbb{R}$ whose transition probability $%
\mathbf{\tilde{P}}$ is specified as follows: for any $\left( \mathbf{x}%
,a\right) \in \mathbb{X}\times \mathbb{R}$ and any bounded Borel function $%
\psi :\mathbb{X}\times \mathbb{R\rightarrow R}$%
\begin{equation*}
\mathbf{\tilde{P}}\psi \left( \mathbf{x},a\right) =\int_{S^{+}}\psi \left(
\mathbf{x}\cdot M,a+\rho (\mathbf{x},M)\right) \mathbf{P}\left( dM\right) .
\end{equation*}%
We denote by $\mathbf{\tilde{P}}_{+}$ the restriction of $\mathbf{\tilde{P}}$
to $\mathbb{X}\times \mathbb{R}_{\ast }^{+}:$%
\begin{equation*}
\mathbf{\tilde{P}}_{+}(\left( \mathbf{x},a\right) ,\cdot )=I\left\{ \cdot
\in \mathbb{X}\times \mathbb{R}_{\ast }^{+}\right\} \mathbf{\tilde{P}}%
(\left( \mathbf{x},a\right) ,\cdot ).
\end{equation*}%
For $\mathbf{x}\in \mathbb{X}$ and $a>0$ denote by $\mathbf{P}_{\mathbf{x}%
,a} $ and $\mathbf{E}_{\mathbf{x},a}$ the probability measure and
expectation on $\left( \Omega ,\mathcal{F},\mathbf{P}\right) $ conditioned
on the event $\left\{ \mathbf{X}_{0}=\mathbf{x},S_{0}=a\right\} $.

Let%
\begin{equation*}
\tau :=\min \left\{ n\geq 1:S_{n}\leq 0\right\} .
\end{equation*}

It is known (see Appendix in \cite{LPP2016}) that if conditions H1-H5 are
valid, then the function $h:$ $\mathbb{X\times R}\rightarrow \mathbb{R}%
_{\ast }^{+}$ defined by
\begin{equation*}
h(\mathbf{x},a):=\lim_{n\rightarrow \infty }\mathbf{E}_{\mathbf{x},a}\left[
S_{n};\tau >n\right]
\end{equation*}%
is $\mathbf{\tilde{P}}_{+}$- harmonic, i.e., for any $\mathbf{x}\in \mathbb{X%
}$ and $a>0$ the equality
\begin{equation}
\mathbf{E}_{\mathbf{x},a}\left[ h(\mathbf{X}_{1},S_{1});\tau >1\right] =h(%
\mathbf{x},a)  \label{HarmonH}
\end{equation}%
is valid.

\begin{lemma}\label{L_abbelow} If conditions H1-H5 are valid, then there exist constants $d>0$ and $C<\infty$ such that for all $\left( \mathbf{x},a\right) \in \mathbb{X}\times \mathbb{R}_{\ast
}^{+}$.
\begin{equation}
  \max\{0,a-d\}<h(\mathbf{x},a)\leq C(1+a)\label{HarmonAbove}
\end{equation}
and
\begin{equation}
  1+a\leq (d+1)(1+h(\mathbf{x},a)).\label{HarmonAbove2}
\end{equation}
\end{lemma}
\textbf{Proof}. For the case when the random walk $\{S_n,\, n\geq 1\}$ is constructed by the iid random matrices taken from $ GL (p,\mathbb{R}),$ the general linear group of $p\times p$ invertible matrices with respect to the ordinary matrix multiplication, and whose distribution satisfies conditions H1-H5, the inequalities (\ref{HarmonAbove}) were established in Theorem 2.2 of \cite{GPP2015}. The proof of (\ref{HarmonAbove}) under our conditions for the semi-group $S^+$ follows the same line if one attracts the arguments included in Appendix of~\cite{LPP2016} (in fact, the upper estimate in (\ref{HarmonAbove}) is estimate (8) in ~\cite{LPP2016}). To check (\ref{HarmonAbove2}) it is sufficient to note that
\begin{equation*}
 a+1=a-d+d+1\leq d+1 + h(\mathbf{x},a)\leq (d+1)(1+h(\mathbf{x},a)).
\end{equation*}
The lemma is proved.

The following important result was proved in \cite{LPP2016}:

\begin{theorem}
\label{T_important}Assume Conditions $H1-H5$. Then, for any $\left( \mathbf{x%
},a\right) \in \mathbb{X}\times \mathbb{R}_{\ast }^{+}$%
\begin{equation*}
\mathbf{P}_{\mathbf{x},a}\left( \tau >n\right) \sim \frac{2}{\sigma \sqrt{%
2\pi n}}\,h(\mathbf{x},a)
\end{equation*}%
as $n\rightarrow \infty $, where $\sigma ^{2}\in (0,\infty )$ is a constant.
Moreover, there exists a constant $c>0$ such that, for any $\left( \mathbf{x}%
,a\right) \in \mathbb{X}\times \mathbb{R}_{\ast }^{+}$%
\begin{equation}
\sqrt{n}\mathbf{P}_{\mathbf{x},a}\left( \tau >n\right) \leq c\left(
1+a\right)  \label{a2.111}
\end{equation}%
for all $n\geq 0$.
\end{theorem}

Let $\mathcal{F}_{n},n\geq 0$ be the $\sigma -$field generated by random
variables
\begin{equation*}
\mathbf{f}_{0},\mathbf{f}_{1},\mathbf{f}_{2},...,\mathbf{f}_{n-1},\mathbf{Z}%
(0),\mathbf{Z}(1),...,\mathbf{Z}(n).
\end{equation*}%
For each $\mathbf{x}\in \mathbb{X}$ and $a>0$ we introduce a new measure $%
\mathbf{\hat{E}}_{\mathbf{x},a}$ on the $\sigma -$field $\mathcal{F}_{\infty
}$ $=\vee _{n=0}^{\infty }\mathcal{F}_{n}$ as follows. For each $n\geq 0$
and each nonnegative random variable $Y_{n}$ measurable with respect to $%
\mathcal{F}_{n}$ we let%
\begin{equation*}
\mathbf{\hat{E}}_{\mathbf{x},a}\left[ Y_{n}\right] :=\frac{1}{h(\mathbf{x},a)%
}\mathbf{E}_{\mathbf{x},a}\left[ Y_{n}h(\mathbf{X}_{n},S_{n});\tau >n\right]
.
\end{equation*}

In view of (\ref{HarmonH}), the new measure $\mathbf{\hat{E}}_{\mathbf{x},a}$
or, what is the same $\mathbf{\hat{P}}_{\mathbf{x},a},$ is well defined on $%
\mathcal{F}_{\infty }$.

The next statement is an analogue of Lemma 2.5 in \cite{4h}.

\begin{lemma}
\label{L_basic}Assume Conditions $H1-H5$. For $k\in \mathbb{N}$, let $Y_{k}$
be a bounded real-valued $\mathcal{F}_{k}$-measurable random variable. Then,
for any $\mathbf{x}\in \mathbb{X}$ and $a>0,$ as $n\rightarrow \infty $%
\begin{equation}
\lim_{n\rightarrow \infty }\mathbf{E}_{\mathbf{x},a}\left[ Y_{k}|\tau >n%
\right] =\mathbf{\hat{E}}_{\mathbf{x},a}\left[ Y_{k}\right] .  \label{FixedK}
\end{equation}%
More generally, let $Y_{1},Y_{2},...$ be a uniformly bounded sequence of
real-valued random variables adopted to the filtration $\mathcal{F}_{\infty}$ which
converges $\mathbf{\hat{P}}_{\mathbf{x},a}$ - a.s. to some random variable $%
Y_{\infty }$ . Then, as $n\rightarrow \infty $%
\begin{equation}
\lim_{n\rightarrow \infty }\mathbf{E}_{\mathbf{x},a}\left[ Y_{n}|\tau >n%
\right] =\mathbf{\hat{E}}_{\mathbf{x},a}\left[ Y_{\infty}\right] .
\label{Growing_n}
\end{equation}
\end{lemma}

\textbf{Proof}. We follow with minor changes the line of proving lemma 2.5
in \cite{4h}. Let%
\begin{equation*}
\mathbf{m}_{x,a}(n):=\mathbf{P}_{\mathbf{x},a}\left( \tau >n\right) .
\end{equation*}%
Clearly,%
\begin{eqnarray*}
\mathbf{E}_{\mathbf{x},a}\left[ Y_{k}|\tau >n\right] &=&\frac{1}{\mathbf{P}_{%
\mathbf{x},a}\left( \tau >n\right) }\mathbf{E}_{\mathbf{x},a}\left[
Y_{k};\tau >n\right] \\
&=&\frac{1}{\mathbf{P}_{\mathbf{x},a}\left( \tau >n\right) }\mathbf{E}_{%
\mathbf{x},a}\left[ Y_{k}\mathbf{m}_{\mathbf{X}_{k},S_{k}}\left( n-k\right)
;\tau >k\right] .
\end{eqnarray*}%
In view of Theorem \ref{T_important}
\begin{equation*}
\lim_{n\rightarrow \infty }\frac{\mathbf{m}_{\mathbf{X}_{k},S_{k}}\left(
n-k\right) }{\mathbf{P}_{\mathbf{x},a}\left( \tau >n\right) }=\frac{h(%
\mathbf{X}_{k},S_{k})}{h(\mathbf{x},a)}\text{ }\quad \mathbf{P}_{\mathbf{x}%
,a}\text{-a.s.}.
\end{equation*}%
Using (\ref{a2.111}) and Lemma \ref{L_abbelow} we see that there exists a constant $C_0>0$ such that%
\begin{equation}
\frac{\mathbf{m}_{\mathbf{X}_{k},S_{k}}\left( n-k\right) }{\mathbf{P}_{%
\mathbf{x},a}\left( \tau >n\right) }\leq C_0\frac{1+S_k}{h(%
\mathbf{x},a)} \leq C_0(1+d)\frac{1+h(\mathbf{X}_{k},S_{k})}{h(%
\mathbf{x},a)}.\label{win}
\end{equation}%
The estimate%
\begin{equation*}
\mathbf{E}_{\mathbf{x},a}\left[ Y_{k}h(\mathbf{X}_{k},S_{k});\tau >k\right]
=h(\mathbf{x},a)\mathbf{\hat{E}}_{\mathbf{x},a}\left[ Y_{k}\right] <\infty
\end{equation*}%
allows us to apply the dominated convergence theorem to get%
\begin{eqnarray*}
\lim_{n\rightarrow \infty }\mathbf{E}_{\mathbf{x},a}\left[ Y_{k}|\tau >n%
\right] &=&\mathbf{E}_{\mathbf{x},a}\left[ Y_{k}\lim_{n\rightarrow \infty }%
\frac{\mathbf{m}_{\mathbf{X}_{k},S_{k}}\left( n-k\right) }{\mathbf{P}_{%
\mathbf{x},a}\left( \tau >n\right) };\tau >k\right] \\
&=&\frac{1}{h(\mathbf{x},a)}\mathbf{E}_{\mathbf{x},a}\left[ Y_{k}h(\mathbf{X}%
_{k},S_{k});\tau >k\right] =\mathbf{\hat{E}}_{\mathbf{x},a}\left[ Y_{k}%
\right] ,
\end{eqnarray*}%
proving (\ref{FixedK}).

To check the validity of (\ref{Growing_n}) fix $\gamma >1$ and observe that
in view of Theorem~\ref{T_important}
\begin{equation*}
\lim_{\gamma \downarrow 1}\lim_{n\rightarrow \infty }\sup \frac{\mathbf{P}_{%
\mathbf{x},a}\left( \tau >n,\tau \leq \gamma n\right) }{\mathbf{P}_{\mathbf{x%
},a}\left( \tau >n\right) }\leq C\lim_{\gamma \downarrow 1}\left( 1-\gamma
^{-1/2}\right) =0
\end{equation*}%
and, by (\ref{win})
\begin{eqnarray*}
\left\vert \mathbf{E}_{\mathbf{x},a}\left[ Y_{n}-Y_{k}|\tau >\gamma n\right]
\right\vert &\leq &\mathbf{E}_{\mathbf{x},a}\left[ \left\vert
Y_{n}-Y_{k}\right\vert \frac{\mathbf{m}_{\mathbf{X}_{n},S_{n}}\left( \left(
\gamma -1\right) n\right) }{\mathbf{P}_{\mathbf{x},a}\left( \tau >\gamma
n\right) };\tau >n\right] \\
&\leq &\frac{C}{h(\mathbf{x},a)}\sqrt{\frac{\gamma }{1-\gamma }}\mathbf{E}_{%
\mathbf{x},a}\left[ \left\vert Y_{n}-Y_{k}\right\vert (1+ h(\mathbf{X}%
_{n},S_{n}));\tau >n\right] \\
&=&C\sqrt{\frac{\gamma }{1-\gamma }}\left(\mathbf{\hat{E}}%
_{\mathbf{x},a}\left[ \left\vert Y_{n}-Y_{k}\right\vert \right]+ \frac{1}{h(\mathbf{x},a)}\mathbf{P}%
_{\mathbf{x},a}\left(\tau>n\right)\right) ,
\end{eqnarray*}%
where $C$ is a positive constant. Letting first $n\rightarrow \infty $ and
then $k\rightarrow \infty $ the right-hand side vanishes by the dominated
convergence theorem. Hence, using the first part of the lemma we see that%
\begin{equation*}
\mathbf{E}_{\mathbf{x},a}\left[ Y_{n};\tau >\gamma n\right] =\left( \mathbf{%
\hat{E}}_{\mathbf{x},a}\left[ Y_{\infty }\right] +o(1)\right) \mathbf{P}_{%
\mathbf{x},a}\left( \tau >\gamma n\right) .
\end{equation*}%
Therefore, for some $C>0$%
\begin{eqnarray*}
&&\lim_{\gamma \downarrow 1}\lim \sup_{n\rightarrow \infty }\frac{\left\vert
\mathbf{E}_{\mathbf{x},a}\left[ Y_{n};\tau >n\right] -\mathbf{\hat{E}}_{%
\mathbf{x},a}\left[ Y_{\infty }\right] \mathbf{P}_{\mathbf{x},a}\left( \tau
>n\right) \right\vert }{\mathbf{P}_{\mathbf{x},a}\left( \tau >n\right) } \\
&&\qquad\leq\lim_{\gamma \downarrow 1}\lim \sup_{n\rightarrow \infty }\frac{%
\left\vert \mathbf{E}_{\mathbf{x},a}\left[ Y_{n};\tau >\gamma n\right] -%
\mathbf{\hat{E}}_{\mathbf{x},a}\left[ Y_{\infty }\right] \mathbf{P}_{\mathbf{%
x},a}\left( \tau >\gamma n\right) \right\vert }{\mathbf{P}_{\mathbf{x}%
,a}\left( \tau >n\right) } \\
&&\quad\qquad+\,C\lim_{\gamma \downarrow 1}\lim \sup_{n\rightarrow \infty }%
\frac{\mathbf{P}_{\mathbf{x},a}\left( \tau >n,\tau \leq \gamma n\right) }{%
\mathbf{P}_{\mathbf{x},a}\left( \tau >n\right) }=0.
\end{eqnarray*}%
The lemma is proved.

Let $\mathbf{f}\left( \mathbf{s}\right) $ be an offspring probability
function. For matrices $A$ and $M$ from $S^{+}$ let
\begin{equation*}
\psi _{\mathbf{f,}A,M}(\mathbf{s}):=\frac{\left\vert A\right\vert }{|A\left(
\mathbf{1}-\mathbf{f}\left( \mathbf{s}\right) \right) |}-\frac{\left\vert
A\right\vert }{|AM\left( \mathbf{1}-\mathbf{s}\right) |}.
\end{equation*}

\begin{lemma}
\label{L_original}If $|B^{(i)}(\mathbf{f})|<\infty $ for all $i\in \left\{
1,...,n\right\} $ and matrix $M=M(\mathbf{f})=\left( m_{kl}\right)
_{k,l=1}^{p}$ satisfies Condition $H3$, then, for any matrix $A=\left(
a_{kl}\right) _{k,l=1}^{p}\in S^{+},\left\vert A\right\vert >0,$ and all $%
\mathbf{s\in }\mathbb{J}\backslash \{\mathbf{1}\}$ the inequality
\begin{equation*}
\psi _{\mathbf{f,}A,M}(\mathbf{s})\leq
bp^{2}\eta \;
\end{equation*}%
is valid.
\end{lemma}

\textbf{Proof}. For $\mathbf{s\in }\mathbb{J}^{p}$ introduce a $p-$%
dimensional vector%
\begin{equation*}
\Delta _{2}\left( \mathbf{s}\right) :=\left( \sum_{j,k}\frac{\partial
^{2}f^{\left( 1\right) }\left( \mathbf{1}\right) }{\partial s_{j}\partial
s_{k}}\left( 1-s_{j}\right) \left( 1-s_{k}\right) ,...,\sum_{j,k}\frac{%
\partial ^{2}f^{\left( p\right) }\left( \mathbf{1}\right) }{\partial
s_{j}\partial s_{k}}\left( 1-s_{j}\right) \left( 1-s_{k}\right) \right) .
\end{equation*}%
Observe that%
\begin{equation}
\left\vert \Delta _{2}\left( \mathbf{s}\right) \right\vert \leq \mu
\left\vert \mathbf{1}-\mathbf{s}\right\vert ^{2}  \label{Norm1}
\end{equation}%
and, in view of Condition $H3$%
\begin{equation*}
\sum_{i=1}^{p}\sum_{l=1}^{p}a_{ik}m_{kl}\geq \min_{k,l}m_{kl}\left\vert
A\right\vert \geq \frac{\left\vert A\right\vert \left\vert M\right\vert }{%
bp^{2}}
\end{equation*}%
implying
\begin{equation}
|AM\left( \mathbf{1}-\mathbf{s}\right) |\geq \frac{\left\vert A\right\vert
\left\vert M\right\vert }{bp^{2}}\left\vert \mathbf{1}-\mathbf{s}\right\vert.
 \label{Norm2}
\end{equation}%
For a fixed $\mathbf{s\in }\mathbb{J}^{p}$ put
\begin{equation*}
H\left( z\right) :=\frac{|A\mathbf{f}\left( \mathbf{s+}z\left( \mathbf{1}-%
\mathbf{s}\right) \right) |}{|A|},\ z\in \left[ 0,1\right] .
\end{equation*}%
Clearly, $H\left( z\right) $ is a probability generating function in $z$ with
\begin{equation*}
H^{\prime }\left( 1\right) =\frac{|AM\left( \mathbf{1}-\mathbf{s}\right) |}{%
|A|},\quad H^{\prime \prime }\left( 1\right) =\frac{|A\Delta _{2}\left(
\mathbf{s}\right) |}{|A|}.
\end{equation*}%
It is known (see, for instance, Lemma 3 in \cite{Kozlov} or Lemma 2.6 in
\cite{GK}), that, for all $z\in \lbrack 0,1)$
\begin{equation*}
0\leq \frac{1}{1-H\left( z\right) }-\frac{1}{H^{\prime }\left( 1\right) (1-z)%
}\leq \frac{H^{\prime \prime }\left( 1\right) }{(H^{\prime }\left( 1\right)
)^{2}}.
\end{equation*}

Further, using simple algebra we see that

\begin{eqnarray*}
0 &\leq &\psi _{\mathbf{f,}A,M}(\mathbf{s})=%
\frac{1}{1-H\left( 0\right) }-\frac{1}{H^{\prime }\left( 1\right) }%
\leq \frac{H^{\prime \prime }\left( 1\right) }{\left(
H^{\prime }\left( 1\right) \right) ^{2}} \\
&=&\frac{|A\Delta _{2}\left( \mathbf{s}\right) |}{|A|}\frac{%
|A|^{2}}{|AM\left( \mathbf{1}-\mathbf{s}\right) |^{2}}=\frac{|A||A\Delta
_{2}\left( \mathbf{s}\right) |}{|AM\left( \mathbf{1}-\mathbf{s}\right) |^{2}}%
.
\end{eqnarray*}%
In view of (\ref{Norm1}) and (\ref{Norm2})
\begin{equation*}
\frac{|A\Delta _{2}\left( \mathbf{s}\right) |}{|AM\left( \mathbf{1}-\mathbf{s%
}\right) |^{2}}\leq \mu \frac{|A|\left\vert \mathbf{1}-\mathbf{s}\right\vert
^{2}bp^{2}}{\left\vert A\right\vert ^{2}\left\vert M\right\vert
^{2}\left\vert \mathbf{1}-\mathbf{s}\right\vert ^{2}}=\mu \frac{bp^{2}}{%
\left\vert A\right\vert \left\vert M\right\vert ^{2}}.
\end{equation*}%
Thus,
\begin{equation*}
0\leq \psi _{\mathbf{f,}A,M}(\mathbf{s})\leq \mu \frac{bp^{2}}{%
\left\vert M\right\vert ^{2}}=bp^{2}\eta,
\end{equation*}%
as desired.

For $0\leq k\leq n$ introduce the notation
\begin{equation*}
\mathbf{f}_{k,n}\left( \mathbf{s}\right) :=\mathbf{f}_{k}(\mathbf{f}%
_{k+1}(...(\mathbf{f}_{n-1}(\mathbf{s}))...))\text{ with }\mathbf{f}%
_{n,n}\left( \mathbf{s}\right) :=\mathbf{s}.
\end{equation*}

\begin{corollary}
\label{C_telescope}For any $\mathbf{x}\in \mathbb{X}$ and $n\geq 0$%
\begin{equation*}
\frac{1}{(\mathbf{x},\mathbf{1}-\mathbf{f}_{0,n}\left( \mathbf{0}\right) )}%
\leq e^{-\ln \left\vert \mathbf{x}R_{n}\right\vert
}+bp^{2}\sum_{k=0}^{n-1}\eta _{k+1}e^{-\ln \left\vert \mathbf{x}%
R_{k}\right\vert }.
\end{equation*}
\end{corollary}

\textbf{Proof}. For $\mathbf{x}=(x_{1},...,x_{p})\in \mathbb{X}$ set $A=A(%
\mathbf{x}):=\left( a_{kl}(\mathbf{x})\right) _{k,l=1}^{p},$ where $a_{kl}(%
\mathbf{x})=x_{k}\delta _{kl}$ and $\delta _{kl}$ is the Kronecker symbol.
Simple algebra gives%
\begin{eqnarray}
\frac{1}{(\mathbf{x},\mathbf{1}-\mathbf{f}_{0,n}\left( \mathbf{s}\right) )}
&=&\frac{1}{|A(\mathbf{1}-\mathbf{f}_{0,n}\left( \mathbf{s}\right) )|}
\notag \\
&=&\frac{1}{|AR_{n}\left( \mathbf{1}-\mathbf{s}\right) |}  \notag \\
&&+\sum_{k=0}^{n-1}\left( \frac{\left\vert AR_{k}\right\vert }{|AR_{k}\left(
\mathbf{1}-\mathbf{f}_{k,n}\left( \mathbf{s}\right) \right) |}-\frac{%
\left\vert AR_{k}\right\vert }{|AR_{k+1}\left( \mathbf{1}-\mathbf{f}%
_{k+1,n}\left( \mathbf{s}\right) \right) |}\right) e^{-\ln \left\vert
AR_{k}\right\vert }  \notag \\
&=&e^{-\ln \left\vert AR_{n}\left( \mathbf{1}-\mathbf{s}\right) \right\vert
}+\sum_{k=0}^{n-1}\psi _{\mathbf{f}_{k}\mathbf{,}AR_{k},M_{k}}(\mathbf{f}%
_{k+1,n}\left( \mathbf{s}\right) )e^{-\ln \left\vert \mathbf{x}%
R_{k}\right\vert }.  \label{Represent}
\end{eqnarray}%
Hence, setting $\mathbf{s}=\mathbf{0}$ yields
\begin{eqnarray}
\frac{1}{(\mathbf{x},\mathbf{1}-\mathbf{f}_{0,n}\left( \mathbf{0}\right) )}
&=&e^{-\ln \left\vert \mathbf{x}R_{n}\right\vert }+\sum_{k=0}^{n-1}\psi _{%
\mathbf{f}_{k}\mathbf{,}AR_{k},M_{k}}(\mathbf{f}_{k+1,n}\left( \mathbf{0}%
\right) )e^{-\ln \left\vert \mathbf{x}R_{k}\right\vert }\quad  \label{FormSurviv}
\\
&\leq&e^{-\ln \left\vert \mathbf{x}R_{n}\right\vert
}+bp^{2}\sum_{k=0}^{n-1}\eta _{k+1}e^{-\ln \left\vert \mathbf{x}%
R_{k}\right\vert },  \notag
\end{eqnarray}%
as it was claimed.

The next statement is an analogue of Lemma 2.7 \ in \cite{4h} and
generalizes to our setting Lemma 3.3 in \cite{LPP2016}.

\begin{lemma}
\label{L_Fourhead1}If the conditions of Theorem \ref{T_main} are valid then,
for any $\mathbf{x}\in \mathbb{X}$ and $a>0$
\begin{equation}
\mathbf{\hat{E}}_{\mathbf{x},a}\left[ \sum_{n=0}^{\infty }\eta
_{n+1}e^{-S_{n}}\right] <\infty .  \label{FiniteSeries}
\end{equation}
\end{lemma}

\textbf{Proof}. We follow the line of proving Lemma 3.3 in \cite{LPP2016}.
First observe that for any $\lambda \in (0,1)$ there exists a constant $%
C=C(\lambda )$ such that $\left( t+1\right) e^{-t}\leq C(\lambda
)e^{-\lambda t}$ for all $t\geq 0$. This estimate, the definition of $%
\mathbf{\hat{E}}_{\mathbf{x},a}$ and (\ref{HarmonAbove}) justify the
inequality%
\begin{eqnarray*}
\mathbf{\hat{E}}_{\mathbf{x},a}\left[ \eta _{n+1}e^{-S_{n}}\right]  &=&\frac{%
1}{h(\mathbf{x},a)}\mathbf{E}_{\mathbf{x},a}\left[ \eta _{n+1}e^{-S_{n}}h(%
\mathbf{X}_{n},S_{n});\tau >n\right]  \\
&\leq &\frac{c}{h(\mathbf{x},a)}\mathbf{E}_{\mathbf{x},a}\left[ \eta
_{n+1}e^{-S_{n}}(1+S_{n});\tau >n\right]  \\
&\leq &\frac{cC(\lambda )}{h(\mathbf{x},a)}\mathbf{E}_{\mathbf{x},a}\left[
\eta _{n+1}e^{-\lambda S_{n}};\tau >n\right] .
\end{eqnarray*}%
Thus,
\begin{eqnarray*}
\mathbf{\hat{E}}_{\mathbf{x},a}\left[ \sum_{n=0}^{\infty }\eta
_{n+1}e^{-S_{n}}\right]  &\leq &\frac{1}{h(\mathbf{x},a)}\mathbf{E}_{\mathbf{%
x},a}\left[ \eta _{1}e^{-a}\right]  \\
&&+\frac{cC(\lambda )}{h(\mathbf{x},a)}\sum_{n=1}^{\infty }\mathbf{E}_{%
\mathbf{x},a}\left[ \eta _{n+1}e^{-\lambda S_{n}};\tau >n\right] .
\end{eqnarray*}%
Further, for any $\kappa ,r\in \mathbb{N}$ with $1\leq r\leq \kappa $ we have%
\begin{eqnarray*}
&&\mathbf{E}_{\mathbf{x},a}\left[ \eta _{n\kappa +r+1}e^{-\lambda S_{n\kappa
+r}};\tau >n\kappa +r\right]  \\
&&\quad \leq \mathbf{E}_{\mathbf{x},a}\left[ e^{-\lambda S_{n\kappa
}};S_{\kappa }>0,S_{2\kappa }>0,...,S_{n\kappa }>0\right] \times \sup_{%
\mathbf{y}\in \mathbb{X}}\mathbf{E}_{\mathbf{y},0}\left[ \eta
_{r+1}e^{-\lambda S_{r}}\right] .
\end{eqnarray*}%
Setting%
\begin{equation*}
\Upsilon _{\kappa }(\mathbf{x},a)=\sum_{n=1}^{\infty }\mathbf{E}_{\mathbf{x}%
,a}\left[ e^{-\lambda S_{n\kappa }};S_{\kappa }>0,S_{2\kappa
}>0,...,S_{n\kappa }>0\right]
\end{equation*}%
we get
\begin{eqnarray*}
&&\sum_{n=1}^{\infty }\mathbf{E}_{\mathbf{x},a}\left[ \eta _{n+1}e^{-\lambda
S_{n}};\tau >n\right]  \\
&\leq &\sum_{r=0}^{\kappa }\mathbf{E}_{\mathbf{x},a}\left[ \eta
_{r+1}e^{-\lambda S_{r}}\right]  \\
&&\quad +\sum_{n=1}^{\infty }\sum_{r=1}^{\kappa }\mathbf{E}_{\mathbf{x},a}%
\left[ \eta _{n\kappa +r+1}e^{-\lambda S_{n\kappa +r}};S_{\kappa
}>0,S_{2\kappa }>0,...,S_{n\kappa }>0\right]  \\
&&\quad \quad \leq \sum_{r=0}^{\kappa }\mathbf{E}_{\mathbf{x},0}\left[ \eta
_{r+1}e^{-\lambda S_{r}}\right]  \\
&&\quad \quad \,+\sum_{n=1}^{\infty }\mathbf{E}_{\mathbf{x},a}\left[
e^{-\lambda S_{n\kappa }};S_{\kappa }>0,S_{2\kappa }>0,...,S_{n\kappa }>0%
\right] \times \sup_{\mathbf{y}\in \mathbb{X}}\sum_{r=1}^{\kappa }\mathbf{E}%
_{\mathbf{y},0}\left[ \eta _{r+1}e^{-\lambda S_{r}}\right]  \\
&&\quad \quad \leq \left( \sup_{\mathbf{y}\in \mathbb{X}}\sum_{r=0}^{\kappa }%
\mathbf{E}_{\mathbf{y},0}\left[ \eta _{r+1}e^{-\lambda S_{r}}\right] \right)
(1+\Upsilon _{\kappa }(\mathbf{x},a)).
\end{eqnarray*}%
Note that according to formula (4.2) in \cite{LPP2016} for any $\delta >0$
there exists a positive integer $\kappa =\kappa (\delta )$ such that%
\begin{equation*}
\mathbf{P}^{\ast \kappa }\left( M:\text{ }\ln \left\vert \mathbf{x}%
M\right\vert \geq \delta \text{ for any }\mathbf{x}\in \mathbb{X}\right) >0.
\end{equation*}%
Moreover, $\Upsilon _{\kappa }(\mathbf{x},a)<\infty $ for the $\kappa $.

Further, setting $\lambda =\varepsilon /\left( 1+\varepsilon \right) ,$
where $\varepsilon $ is the same as in (\ref{SecondFinite}), we deduce by H%
\"{o}lder inequality that
\begin{eqnarray*}
\mathbf{E}_{\mathbf{y},0}\left[ \eta _{r+1}e^{-\varepsilon
S_{r}/(1+\varepsilon )}\right] &\leq &\left( \mathbf{E}_{\mathbf{y},0}\left[
 \eta _{r+1}^{1+\varepsilon }\right] \right)
^{1/(1+\varepsilon )}\left( \mathbf{E}_{\mathbf{y},0}\left[ e^{-S_{r}}\right]
\right) ^{\varepsilon /(1+\varepsilon )} \\
&\leq &\left( \mathbf{E}\left[ \eta _{r+1}^{1+\varepsilon }%
\right] \right) ^{1/(1+\varepsilon )}\left( \sup_{\mathbf{y}\in \mathbb{X}}%
\mathbf{E}_{\mathbf{y},0}\left[ e^{-S_{1}}\right] \right) ^{r\varepsilon
/(1+\varepsilon )}
\end{eqnarray*}%
and the right-hand side is finite by conditions (\ref{ExponFinite}) and (\ref%
{SecondFinite}).

Combining the estimates above gives (\ref{FiniteSeries}).

 Having Lemma \ref{L_Fourhead1} in hands and using basically the arguments of \cite{GK} and \cite{LPP2016} we are now ready to prove the main result of the note.

\textbf{Proof of Theorem \ref{T_main}}. For a fixed $i\in \left\{
1,...,p\right\} $ we write
\begin{eqnarray*}
\mathbf{P}\left( \mathbf{Z}(n)\neq \mathbf{0}|\mathbf{Z}(0)=\mathbf{e}%
_{i}\right) &=&\mathbf{E}\left[ 1-f_{0,n}^{(i)}\left( \mathbf{0}\right) %
\right] \\
&=&\mathbf{E}\left[ \left( \mathbf{e}_{i},\mathbf{1}-\mathbf{f}_{0,n}\left(
\mathbf{0}\right) \right) \right] \\
&=&\mathbf{E}_{\mathbf{e}_{i},a}\left[ \left( \mathbf{e}_{i},\mathbf{1}-%
\mathbf{f}_{0,n}\left( \mathbf{0}\right) \right) ;\tau >n\right] \\
&&+\mathbf{E}_{\mathbf{e}_{i},a}\left[ \left( \mathbf{e}_{i},\mathbf{1}-%
\mathbf{f}_{0,n}\left( \mathbf{0}\right) \right) ;\tau \leq n\right] .
\end{eqnarray*}

In view of the inequalities
\begin{equation*}
\mathbf{1}-\mathbf{f}_{0,n}\left( \mathbf{0}\right) \leq \mathbf{1}-\mathbf{f%
}_{0,k}\left( \mathbf{0}\right) \leq M_{0}M_{1}\cdot \cdot \cdot M_{k-1}%
\mathbf{1}=R_{k}\mathbf{1}\text{,}
\end{equation*}%
being valid for all $0\leq k\leq n$ we have for a constant $C_{1}\in
(0,\infty ):$
\begin{multline*}
\mathbf{E}_{\mathbf{e}_{i},a}\left[ \left( \mathbf{e}_{i},\mathbf{1}-\mathbf{%
f}_{0,n}\left( \mathbf{0}\right) \right) ;\tau \leq n\right] \leq C_{1}%
\mathbf{E}_{\mathbf{e}_{i},a}\left[ \min_{0\leq k\leq n}|\mathbf{e}_{i}R_{k}%
\mathbf{1|};\tau \leq n\right] \\
=C_{1}\mathbf{E}_{\mathbf{e}_{i},0}\left[  \min_{0\leq k\leq n}|\mathbf{e}_{i}R_{k}%
\mathbf{1|};\min_{0\leq k\leq n}\ln
\left\vert \mathbf{e}_{i}R_{k}\right\vert \leq -a\right] .
\end{multline*}%
To evaluate the right-hand side we use the inequalities%
\begin{multline*}
\mathbf{E}_{\mathbf{e}_{i},0}\left[e^{\min_{0\leq k\leq n}\ln|\mathbf{e}_{i}R_{k}%
\mathbf{1|}};\min_{0\leq k\leq n}\ln \left\vert \mathbf{e%
}_{i}R_{k}\right\vert \leq -a\right]\\
 \leq \sum_{j\geq a}^{\infty }e^{-j}%
\mathbf{P}_{\mathbf{e}_{i},0}\left( -j<\min_{0\leq k\leq n}\ln \left\vert
\mathbf{e}_{i}R_{k}\right\vert \leq -j+1\right) \\
\leq \sum_{j\geq a}^{\infty }e^{-j}\mathbf{P}_{\mathbf{e}_{i},j}\left( \tau
>n\right)
\leq \frac{C}{\sqrt{n}}\sum_{j\geq a}^{\infty }e^{-j}(j+1),
\end{multline*}%
where the last estimate follows from (\ref{a2.111}).

Note that the estimates above imply
\begin{equation*}
\lim_{a\rightarrow \infty }\lim \sup_{n\rightarrow \infty }\sqrt{n}\mathbf{E}%
_{\mathbf{e}_{i},a}\left[ \left( \mathbf{e}_{i},\mathbf{1}-\mathbf{f}%
_{0,n}\left( \mathbf{0}\right) \right) ;\tau \leq n\right] =0.
\end{equation*}

Set $\mathbf{f}_{k,\infty }\left( \mathbf{0}\right) :=\lim_{n\rightarrow
\infty }\mathbf{f}_{k,n}\left( \mathbf{0}\right) $. It is not difficult to
deduce from Lemma \ref{L_Fourhead1} that $\ln \left\vert \mathbf{e}%
_{i}R_{n}\right\vert \rightarrow \infty $ \ $\mathbf{\hat{P}}_{\mathbf{e}%
_{i},a}-$a.s. for any $a>0$. Recalling formula (\ref{FormSurviv}) we see
that, as $n\rightarrow \infty $%
\begin{eqnarray*}
\frac{1}{\left( \mathbf{e}_{i},\mathbf{1}-\mathbf{f}_{0,n}\left( \mathbf{0}%
\right) \right) } &\rightarrow &\frac{1}{\left( \mathbf{e}_{i},\mathbf{1}-%
\mathbf{f}_{0,\infty }\left( \mathbf{0}\right) \right) } \\
&=&\sum_{k=0}^{\infty }\psi _{\mathbf{f}_{k}\mathbf{,}AR_{k},M_{k}}(\mathbf{f%
}_{k+1,\infty }\left( \mathbf{0}\right) )e^{-\ln \left\vert \mathbf{e}%
_{i}R_{k}\right\vert }
\end{eqnarray*}%
$\mathbf{\hat{P}}_{\mathbf{e}_{i},a}-$a.s., where now $A=A(\mathbf{e}%
_{i})=\left( a_{kl}\right) _{k,l=1}^{p}$ is a $p\times p$ matrix with $%
a_{ii}=1$ and all other elements equal to zero. Thus, as $n\rightarrow
\infty $
\begin{equation*}
Y_{n}:=\left( \mathbf{e}_{i},\mathbf{1}-\mathbf{f}_{0,n}\left( \mathbf{0}%
\right) \right) \rightarrow \left( \mathbf{e}_{i},\mathbf{1}-\mathbf{f}%
_{0,\infty }\left( \mathbf{0}\right) \right) \text{ }:=Y_{\infty }\text{\ \ }%
\mathbf{\hat{P}}_{\mathbf{e}_{i},a}-\text{a.s.}
\end{equation*}%
and, by Corollary \ref{C_telescope} and Lemma \ref{L_Fourhead1}
\begin{equation*}
Y_{\infty }\geq \left( \sum_{k=0}^{\infty }\eta _{k+1}e^{-\ln \left\vert
\mathbf{e}_{i}R_{k}\right\vert }\right) ^{-1}>0
\end{equation*}%
\ $\mathbf{\hat{P}}_{\mathbf{e}_{i},a}-$a.s.$.$ Since $0\leq Y_{n}\leq 1$,
applying Lemma \ref{L_basic} gives%
\begin{equation}
\lim_{n\rightarrow \infty }\sqrt{n}\mathbf{E}_{\mathbf{e_{i}},a}\left[
Y_{n};\tau >n\right] =\frac{2}{\sigma \sqrt{2\pi }}\mathbf{\hat{E}}_{\mathbf{%
e}_{i},a}\left[ Y_{\infty }\right] h(\mathbf{e}_{i},a)>0.  \label{FirstStep}
\end{equation}%
The left-hand side of (\ref{FirstStep}) is increasing in $a$. Hence, the
right-hand side of (\ref{FirstStep}) is also increasing in $a$. Therefore,
the limit
\begin{equation*}
\beta _{i}:=\lim_{a\rightarrow \infty }\frac{2}{\sigma \sqrt{2\pi }}\mathbf{%
\hat{E}}_{\mathbf{e}_{i},a}\left[ Y_{\infty }\right] h(\mathbf{e}_{i},a)
\end{equation*}%
exists. It is finite in view of (\ref{FirstStep}) and the estimate
\begin{equation*}
\mathbf{E}_{\mathbf{e}_{i},a}\left[ \left( \mathbf{e}_{i},\mathbf{1}-\mathbf{%
f}_{0,n}\left( \mathbf{0}\right) \right) ;\tau \leq n\right] \leq \frac{%
CC_{1}}{\sqrt{n}}\sum_{k\geq a}^{\infty }e^{-k}(k+1).
\end{equation*}%
Thus,
\begin{equation*}
\lim_{n\rightarrow \infty }\sqrt{n}\mathbf{P}\left( \mathbf{Z}(n)\neq
\mathbf{0}|\mathbf{Z}(0)=\mathbf{e}_{i}\right) =\beta _{i}.
\end{equation*}

Theorem \ref{T_main} is proved.

\end{document}